# AN ASTRONOMICAL VIEW-POINT ON THE EASTER DATE

Ovidiu Vaduvescu

PhD Student, York University, Toronto ON Canada

Research Associate, The Astronomical Institute, Bucharest Romania

ABSTRACT

Old algorithms are still in use nowadays setting the Easter date accordingly to some elements giving the Spring Equinox and the New Moon for a given year. At least two different approaches (by the Orthodox and Catholic Churches) provide different results, most of the times different than the actual moments given by the two astronomical phenomena. A new algorithm built accordingly to the Niceea rule and based on some simple astronomical formulas is presented. It is proposed to replace the old different approaches, in order to celebrate the Easter at the same common date. A Windows program is provided, together with a table listing the Easter date between the years 1950-2050 calculated using the Orthodox, Catholic and the Astronomical approaches.

1. INTRODUCTION

The rule promulgated in Niceea (Bitinia, Minor Asia) in the year 325 A.D. concerning the Easter date sets this holiday "*On the Sunday which follows the 14-th day of the Moon which reaches this age on 21 March or immediately after that*" (e.g. Levy, 1974). In this definition, by "the 14-th day of the Moon" one has to understand "the age 14 of the Moon" (i.e. the day 14 days after the New Moon).

It is very difficult to establish the origin of the first algorithms used to calculate the Easter date. It is supposed that The Holy Father Chiril, bishop of Alexandria, was the first who approached this problem between 430-440 A.D. The first systematic algorithm dates from the year 456 A.D. and might be assigned to Victorius d'Aquitania.

The elementary arithmetic (involving only summation, multiplication and division) represents the instrument to calculate the Easter date accordingly to the classical approach presented in many papers (e.g. Levy 1975). The elements required by this algorithm (Vaduvescu, 1993, Vaduvescu and Soulie, 1994) are expressed as "modulo" functions of some arithmetical progression terms.

The imperfection of this method consists in the fact that some entire numbers can not represent accurately, for long time periods, the natural phenomena which set the Easter date, namely the Vernal Equinox and the New Moon phase. Moreover, because of the Reformation of the Calendar in 1582, the elements of the algorithm were modified. While the Catholic Church approach defines these elements based on the Gregorian calendar, the Orthodox Church preserves in its algorithm the elements expressed in the old Julian calendar, correcting the result by a number of days equal to the difference between the two calendars (Vaduvescu and Soulie, 1994). Nevertheless, although more precise, the Gregorian calendar will generate in the year 4317 a delay of a day in comparison with the tropical time (BdL, 1970, 1971).

This is my third and probably the last paper dealing with the Easter date. It was originally written as a manuscript in 1994. While the second paper (Vaduvescu and Soulie, 1994) provides an algorithm to calculate the Easter using elements defined based on the Gregorian Calendar, the first paper (Vaduvescu, 1993) introduces the computational elements necessary to calculate the Easter by the Orthodox Church, it compares them with the Catholic Church approach, providing a program to calculate both. In the present paper we propose a new algorithm based on some classic astronomical formulae giving the Equinox date and the Moon phase, instead of using the old elements set to different calendars to calculate the Easter Date.

2. THE VERNAL EQUINOX

It is well known that setting the equinox each year on 21 March represents an approximation, this phenomenon varying in time, as a function of the moment when the longitude of the Sun is zero. Savoie (1988) showed that in the year 325 A.D. the equinox has fallen on 20 March, about around midday (UT). Therefore, an astronomical view-point on the rigorous treating of the problem would be to define the Easter date as a function of the equinox, and not of 21 March.

Accordingly with Newcomb's solar theory (e.g. Meeus, 1984, Meeus, 1991), denoting by *JJ* the Julian date, let *T* be the following fraction:

$$T = (JJ - 2415020.0) / 36525$$

Then, the mean longitude of the Sun is given by the following truncated series:

$$L = 279°.69668 + 36000°.76892\, T + 0°.0003025\, T^2 \qquad (1)$$

The mean anomaly of the Sun is:

$$M = 358°.47583 + 35999°.04975\,T - 0°,000150\,T^2 - 0°.0000033\,T^3 \qquad (2)$$

and the equation of the Sun center is:

$$C = (1°.919460 - 0°.004789\,T - 0°.000014\,T^2)\sin(M) +$$
$$+ (0°.020094 - 0°.0001T)\sin(2M) + 0°.000293\sin(3M) \qquad (3)$$

Then, the true longitude of the Sun is:

$$TL = L + C \qquad (4)$$

Because TL is a continuous and monotonous function of the JJ, the method of the bisection or Newton's method for example can provide the moment denoted by $JJ_e$ of the year when the mean longitude of the Sun is a multiple of 360°, considering for example an initial value $JJ_0$ corresponding to 21 March 0h UT.

3. COMMON YEAR OR LEAP YEAR

Let *m* be the year for which we search the Easter date. The number of the days in February, expressed in the Gregorian calendar, is 28+nf, where we denote by *nf* the following number:

$$\begin{aligned}&\text{if } [m]_{400} = 0, \text{ then } nf = 1\\&\text{if } [m]_{400} <> 0 \text{ and } [m]_{100} = 0, \text{ then } nf = 0\\&\text{if } [m]_{100} <> 0 \text{ and } [m]_4 = 0, \text{ then } nf = 1 \qquad (5)\\&\text{if } [m]_{100} <> 0 \text{ and } [m]_4 <> 0, \text{ then } nf = 0\end{aligned}$$

where $[x]_y$ represents the rest of the entire division of the entire numbers *x* by *y*.

4. THE NEW MOON

The fraction of year corresponding to 1 March is:

$$f = (59+nf)/(365+nf) \qquad (6)$$

by means of which the multiplication constant corresponding to the mean phase of New Moon is:

$$k = \text{int}[(m + f - 1900.0) \times 12.3685] + 1 \qquad (7)$$

where by "int[x]" be denoted the entire part of the real *x*.

According to Meeus (1986), the moment of mean phase of the New Moon is:

$$JJ_{NM} = 2415020.75933 + 29.53058868\ k \qquad (8)$$

which sets the moment of the first New Moon subsequent to 1 March.

5. ADDITIONAL ACCURACY

Formula (1) can be improved by taking into consideration planetary perturbations, while (8) by taking into account the aberration of the Sun and the real phases, according to Meeus (1986, 1991). Also, corrections of diurnal parallax and transformations of local time corresponding to the coordinates of Niceea could be added if necessary to the study. We do not insist on these.

6. THE 14-TH DAY OF THE MOON

If $(JJ_{NM}+14) >= JJ_e$, then this is the Paschal Moon, otherwise we are placed into the former cycle. So, in the second case, by assigning:

$$JJ_{NM} = JJ_{NM} + 29.53058868$$

we will be placed in the Paschal Moon.

Denoting by $JJ_{1J}$ the Julian date corresponding to 1 January, one calculates the expression:

$$z_a = \text{int}[JJ_{NM} - JJ_{1J}] + 1 \qquad (9)$$

which represents the day of the year corresponding to the Paschal New Moon.

7. THE NEXT SUNDAY

The first Sunday after the $j$-th day of an year $m$ is given by:

$$d = j + [2 - j - m_{ac}]_7 \qquad (10)$$

where by "$m_{ac}$" we denoted "the Catholic hand of the year" (the Catholic element equivalent to the Orthodox dominical letter element) given by:

$$m_{ac} = [u + \text{int}(u/4) + \text{int}(c/4) - 2c]_7 + 1 \qquad (11)$$

where:

$$c = \text{int}(m/100) \text{ and } u = [m]_{100}$$

If $m$ is a leap year (nf=1), then $m_{ac}$ is obtained by assigning to $m_{ac} = (m_{ac} - 1)$.

Finally, the day of the year corresponding to the Easter is:

$$z_P = z_a + 14 + [2 - z_a - m_{ac}]_7 \qquad (12)$$

which remains to be transformed into calendar date:

if $z_P > (90+nf)$, then the Easter falls on $(z_P - 90 - nf)$ April
if $z_P <= (90+nf)$, then the Easter falls on $(z_P - 59 - nf)$ March $\qquad (13)$

## 8. EXAMPLE

Let m=1994. One gets $JJ_0 = 2449432.5$; $JJ_{1J} = 2449353.5$ and with (1), (2), (3), (4), by means of the half interval method with $JJ_1=2449432.0$ after four recurrences, one gets TL=0.02381 and $JJ_e=2449432.375$ (that is 20 March 21h UT, in comparison with 20 March 20.5h UT from the ephemeris). (5) yields $nf = 0$. The new moon is given by (6), (7) and (8): $f = 0.1616438$; $k = 1165$; $JJ_{NM} = 2449423.895$ (that is 12 March 9h UT, in comparison with 12 March 7h UT from the ephemeris). Because $(JJ_{NM}+14) >= JJ_e$, we are placed in the searched cycle, then from (9) $z_a=71$. From (11), $m_{ac} = 7$, then accordingly with (12), $z_P = 86$. Consequently, the astronomical Easter falls on 27 March 1994, in comparison with 3 April the Catholic Easter or 1 May the Orthodox Easter.

## 9. PROGRAM, COMPARISON

We wrote a PC code in Delphi 4 for Windows based on the above algorithm which provides the Easter date between two limits, given by a starting year and a number of years. The executable is available online (Vaduvescu, 2004).

A table giving the Astronomical Easter between the years 1950-2050 (calculated with the above program), the Orthodox Easter, and the Catholic Easter (both provided by the program in Vaduvescu, 1993) is shown in Table 1 in the Appendix.

Major differences can be seen from Table 1 in most cases between the Astronomical Easter date (supposed to be calculated using the most correct algorithm) and the Orthodox approach (which has at its base the Julian calendar).

From the 101 years analyzed, there are 36 cases (35.6%) for which some differences can be seen between the Astronomical Easter date and the Catholic algorithm (which employs the Gregorian calendar to calculate the elements used for Easter algorithm). No apparent periodicity can be seen. Interestingly, these differences represent only two intervals (in the sense Astronomic minus Catholic): minus one week, appearing in 21 cases (20.8% from the total 101 years), and plus four weeks (in the same sense) appearing in the rest 15 cases (14.8% cases).

APPENDIX

Table 1: Easter dates calculated between 1950-2050 (first column)
based on the Orthodox (column 2), Catholic (column 3), and the Astronomical approach.

```
YEAR      ORTHODOX     CATHOLIC     ASTRONOMIC
==========================================
1950       9 April      9 April      2 April
1951      29 April     25 March     22 April
1952      20 April     13 April     13 April
1953       5 April      5 April     29 March
1954      25 April     18 April     18 April
1955      17 April     10 April     10 April
1956       6 May        1 April      1 April
1957      21 April     21 April     14 April
1958      13 April      6 April      6 April
1959       3 May       29 March     26 April
1960      17 April     17 April     10 April
1961       9 April      2 April      2 April
1962      29 April     22 April     22 April
1963      14 April     14 April     14 April
1964       3 May       29 March     29 March
1965      25 April     18 April     18 April
1966      10 April     10 April     10 April
1967      30 April     26 March     26 March
1968      21 April     14 April     14 April
1969      13 April      6 April      6 April
1970      26 April     29 March     26 April
1971      18 April     11 April     11 April
1972       9 April      2 April      2 April
1973      29 April     22 April     22 April
1974      14 April     14 April      7 April
1975       4 May       30 March     30 March
1976      25 April     18 April     18 April
1977      10 April     10 April      3 April
1978      30 April     26 March     23 April
1979      22 April     15 April     15 April
1980       6 April      6 April     30 March
1981      26 April     19 April     19 April
1982      18 April     11 April     11 April
1983       8 May        3 April      3 April
1984      22 April     22 April     15 April
1985      14 April      7 April      7 April
1986       4 May       30 March     27 April
1987      19 April     19 April     12 April
1988      10 April      3 April      3 April
1989      30 April     26 March     23 April
1990      15 April     15 April     15 April
1991       7 April     31 March     31 March
1992      26 April     19 April     19 April
1993      18 April     11 April     11 April
1994       1 May        3 April     27 March
1995      23 April     16 April     16 April
1996      14 April      7 April      7 April
1997      27 April     30 March     27 April
1998      19 April     12 April     12 April
1999      11 April      4 April      4 April
2000      30 April     23 April     23 April
2001      15 April     15 April      8 April
```

Table 1 (continued): Easter dates calculated between 1950-2050 (first column) based on the Orthodox (column 2), Catholic (column 3), and the Astronomical approach.

```
YEAR       ORTHODOX     CATHOLIC     ASTRONOMIC
==========================================
2002        5 May       31 March     31 March
2003       27 April     20 April     20 April
2004       11 April     11 April      4 April
2005        1 May       27 March     24 April
2006       23 April     16 April     16 April
2007        8 April      8 April      1 April
2008       27 April     23 March     20 April
2009       19 April     12 April     12 April
2010        4 April      4 April      4 April
2011       24 April     24 April     17 April
2012       15 April      8 April      8 April
2013        5 May       31 March     31 March
2014       20 April     20 April     13 April
2015       12 April      5 April      5 April
2016        1 May       27 March     24 April
2017       16 April     16 April     16 April
2018        8 April      1 April      1 April
2019       28 April     21 April     21 April
2020       19 April     12 April     12 April
2021        2 May        4 April     28 March
2022       24 April     17 April     17 April
2023       16 April      9 April      9 April
2024        5 May       31 March     28 April
2025       20 April     20 April     13 April
2026       12 April      5 April      5 April
2027        2 May       28 March     25 April
2028       16 April     16 April      9 April
2029        8 April      1 April      1 April
2030       28 April     21 April     21 April
2031       13 April     13 April      6 April
2032        2 May       28 March     28 March
2033       24 April     17 April     17 April
2034        9 April      9 April      9 April
2035       29 April     25 March     22 April
2036       20 April     13 April     13 April
2037        5 April      5 April      5 April
2038       25 April     25 April     18 April
2039       17 April     10 April     10 April
2040        6 May        1 April      1 April
2041       21 April     21 April     21 April
2042       13 April      6 April      6 April
2043        3 May       29 March     26 April
2044       24 April     17 April     17 April
2045        9 April      9 April      2 April
2046       29 April     25 March     22 April
2047       21 April     14 April     14 April
2048        5 April      5 April     29 March
2049       25 April     18 April     18 April
2050       17 April     10 April     10 April
```